\documentclass[12pt]{conm-p-l} 
\usepackage{geometry}
\usepackage{graphicx} 
\usepackage{amssymb}
\usepackage{epstopdf}
\usepackage{amsmath,amscd}
\usepackage{amsthm} 
\usepackage{url,verbatim}
\usepackage{mathtools}

\RequirePackage[colorlinks,citecolor=blue,urlcolor=blue]{hyperref}
\usepackage{breakurl} 
\theoremstyle{plain}
\newtheorem{theorem}{Theorem}
 
\newtheorem{lemma}[theorem]{Lemma}

\newtheorem{remark}[theorem]{Remark}

\newcounter{mycount}

\newenvironment{letlist}{\begin{list}{\rm(\alph{mycount})}%
{\usecounter{mycount}\labelwidth=1cm\itemsep 0pt}}{\end{list}}

\numberwithin{equation}{section} 
\numberwithin{theorem}{section}
\numberwithin{figure}{section}

\newcommand\HH{{\mathbb H}}

\newcommand\qq{\qquad} 
\newcommand\q{\quad} 
\newcommand\si{\sigma}

\newcommand\Si{\Sigma}
\newcommand\g{\gamma}

\newcommand\De{\Delta}

\newcommand\CC{{\mathbb C}}

\newcommand\ZZ{{\mathbb Z}} 
\newcommand\VV{{\mathbb V}}
\newcommand\EE{{\mathbb E}}

\newcommand\La{\Lambda} 

\newcommand\eps{\epsilon} 
\newcommand\ot{1-2\ } 
\newcommand\resp{respectively}
 
\newcommand\lra{\leftrightarrow}
\newcommand\oo{\infty} 
 
\newcommand\TT{{\mathbb T}}
 
\newcommand\de{\delta} 

\newcommand\Pf{\mathrm{Pf}\, } 

\newcommand\rv{\mathrm{v}} 
\newcommand\re{\mathrm{e}}

\newcommand\HnD{\HH_{n,\De}}

\newcommand\mo{\mu_\oo} 
\newcommand\Pip{\Pi^{\text{\rm poly}}}
\newcommand\Znp{Z_n^{\text{\rm poly}}}
\newcommand\Pipi{\Pi_{e,f}}
\newcommand\ac{a_{\mathrm{c}}}
\newcommand\pa{\ell}

\title[The \ot model]{The \ot model}

\author{Geoffrey R.\ Grimmett} \address{Statistical Laboratory, Centre for
Mathematical Sciences, Cambridge University, Wilberforce Road, Cambridge CB3
0WB, UK} 
\email{g.r.grimmett@statslab.cam.ac.uk}
\urladdr{\url{http://www.statslab.cam.ac.uk/~grg/}}

\author{Zhongyang Li}
\address{Department of Mathematics,
University of Connecticut, Storrs, Connecticut 06269-3009, USA}
\email{zhongyang.li@uconn.edu}
\urladdr{\url{http://www.math.uconn.edu/~zhongyang/}}

\begin{document}

\begin{abstract} 
The current paper is a short review of rigorous results for
the \ot model.
The \ot model on the hexagonal lattice is a model of
statistical mechanics in which each vertex is constrained to have degree either
$1$ or $2$. It was proposed in a study by Schwartz and Bruck of constrained coding
systems, and is strongly connected to the dimer model on a decoration of the lattice,
and to an enhanced Ising model and an associated polygon model on
the graph derived from the hexagonal lattice by adding a further vertex in the middle of each edge.
 
The general \ot model possesses three parameters
$a$, $b$, $c$. The fundamental technique is to represent
probabilities of interest as ratios of counts of dimer coverings of certain associated 
graphs, and to apply the Pfaffian method of Kasteleyn, Fisher, and Temperley.

Of special interest is the existence (or not) 
of phase transitions. It turns out that all clusters of
the infinite-volume limit are almost surely finite.
On the other hand, the existence (with strictly positive probability)
of infinite `homogeneous' clusters, containing vertices of given type,
depends on the values of the parameters.

A further type of phase transition emerges in the study of the
two-edge correlation function, and in this case the critical surface may be found explicitly.
For instance, when $a \ge b \ge c > 0$, the surface given by
$\sqrt a = \sqrt b + \sqrt c$ is critical. 
\end{abstract}

\date{July 14, 2015, revised September 7, 2015} 
\keywords{\ot model, dimers, polygon model, Ising model, perfect
matching, Kasteleyn matrix, phase transition, percolation.}
\subjclass[2010]{82B20, 60K35, 05C70}

\maketitle

\section{Origin of the \ot model}\label{12d}

The \ot model originated in the work of computer scientists Schwartz and Bruck \cite{SB08}
on constrained coding systems. They studied  an array of variables 
on the hexagonal lattice $\HH$ subject to the \lq not all equal' constraint. Of particular
interest to them was the asymptotic behaviour of the number of acceptable configurations
on  large bounded regions $\HH_n$ of the lattice, in the `thermodynamic limit'
as $\HH_n\uparrow \HH$.
Using the method of so-called `holographic reduction', they were able
to map their counting problem to one of counting the number of perfect
matchings (or `dimer coverings') on a certain graph derived from the hexagonal lattice. This last problem
may be solved using the Pfaffian representation of Kasteleyn \cite{Kast61}, 
Fisher \cite{F61}, and Temperley and Fisher \cite{TF61}.

\begin{figure}[htbp]
\centering
\scalebox{1}[1]{\includegraphics*[width=0.6\hsize]{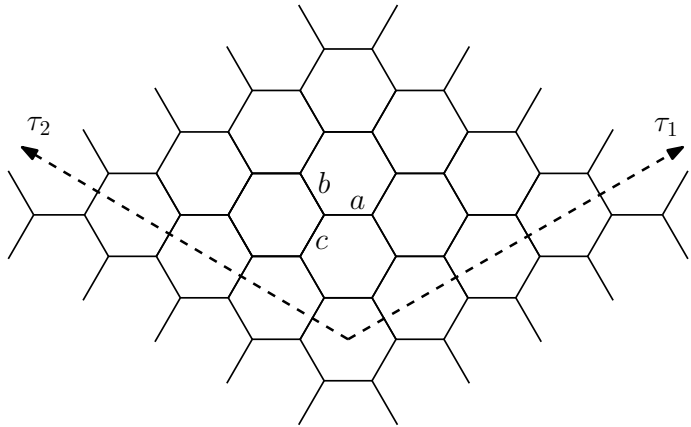}}
\caption{The graph $\HH_n$ is an $n\times n$ lozenge wrapped onto a torus.
A horizontal (\resp, NW/SE, NE/SW) edge is said to be of type
$a$ (\resp, type $b$, $c$).}\label{fig:hex0}
\end{figure}

When rephrased in the language of statistical mechanics, the work of Schwartz and Bruck
amounts to the calculation of the partition function of the following \emph{\ot model}
of probability theory and mathematical physics.
Let $\HH=(\VV,\EE)$ be the hexagonal lattice of Figure \ref{fig:hex0}, and
let $\Si=\{-1,1\}^\EE$ be the set of configurations of absent/present edges,
where the local state $-1$ (\resp, $1$) means \emph{absent} (\resp, \emph{present}).
The sample space is the subset of $\Si$ containing
all $\si\in\Si$ such that: every vertex of $\HH$ is incident to either
one or two present edges. Thus, a configuration comprises disjoint paths and cycles
of present edges.

Let $\HH_n$ be an $n\times n$ subgraph of $\HH$
with periodic boundary conditions, and let $\mu_n$ be the uniform probability measure
on  the set of \ot configurations on $\HH_n$.
We ask for properties of $\mu_n$ in the limit as $n\to\oo$.
In particular, does the limit measure $\mo:=\lim_{n\to\oo}\mu_n$ exist, and, if so,
what can be said about the long-range correlations of edge-states under $\mo$?
It turns out that the connection to dimers may be exploited to answer such questions.

The above system is a lattice model
whose partition function can be computed by the calculation of certain determinants using the
holographic algorithm of Valiant \cite{Val}. 
By introducing an invertible $2\times 2$ matrix on edges of a graph and
conducting a base change,  the partition function of a general vertex
model on a graph $G$ is transformed  into the partition function of perfect matchings on a certain decorated
version of $G$. 
Valiant's original holographic algorithm can be generalized by
assigning different bases to different edges (see \cite{ZLvm}), and the ensuing 
algorithm can be used to compute partition functions of a larger class of
vertex models in polynomial time.

The holographic algorithm seems, however, not to be the most efficient way to solve the \ot model. 
In particular, the correspondence between the
\ot partition function and the dimer
partition function on a corresponding Fisher graph, via the
base change method, is not measure-preserving; thus, the computation of local
statistics and related probabilities becomes complicated, even if possible.
An alternative measure-preserving correspondence was introduced in
\cite{ZL2}, and this permits a number of representations in closed form
of probabilities associated with the \ot model. This method,
and some of its consequences, will be described in the current review.

Certain properties of the underlying hexagonal lattice are utilized heavily in this work, such as
trivalence, planarity, and support of a $\ZZ^2$ action. It may be possible
to extend the results summarized here to certain other graphs with such properties, including
the Archimedean $(3, 12^2)$ and $(4, 8^2)$ lattices.

The formal definition of the \ot model is presented in Section \ref{sec:def}. The model
has strong connections to the dimer and Ising models as well as to a certain polygon 
model, and these connections are laid out in Section \ref{sec:conn}.  Two approaches 
to the issue of phase transition are outlined in Section \ref{sec:pt}, using  
the geometry and the correlation structure, \resp, and an exact formula for the
critical surface is given in the second case.

\section{Definition of the \ot model}\label{sec:def}

Whereas the \ot model of \cite{SB08} is uniform in that there is only one parameter, 
we present here the more general three-parameter model of \cite{ZLvm}.

Let $n \ge 1$, and let $\tau_1$, $\tau_2$ be the two shifts of $\HH$ as in Figure \ref{fig:hex0}.
The pair $(\tau_1,\tau_2)$ generates a $\ZZ^2$ action on $\HH$, and we write 
$\HH_n=(V_n,E_n)$ for
the (toroidal) quotient graph of $\HH$ under the subgroup of $\ZZ^2$ generated 
by the powers $\tau_1^n$ and $\tau_2^n$.  The 
configuration space $\Si_n$ is the set of all $\si\in\{-1,1\}^{E_n}$ such
that every $v \in V_n$ is incident to either $1$ or $2$ edges $e$ with $\si(e)=1$.
Note that $\si\in\Si_n$ if and only if $-\si\in\Si_n$. It is sometimes
convenient to work with the vector $\si'$ given by $\si'(e)=\frac12(1+\si(e))$.

\begin{figure}[htbp] 
\centerline{\includegraphics*[width=0.98\hsize]{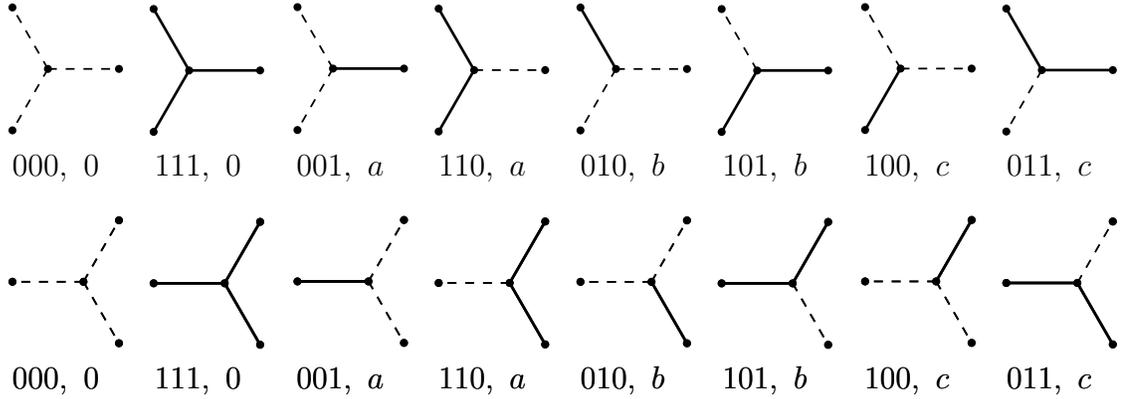}}
\caption{Each vertex has eight possible local configurations,
six of which have nontrivial weights, $a$, $b$, or $c$ as above.
A vertex $v$ is said to be of \emph{type} $s$ in configuration $\si$
if $\si|_v$ has weight $s\in\{a,b,c\}$.
Whereas the type of an edge is \emph{deterministic} (see Figure \ref{fig:hex0}),
the type of a vertex is \emph{random}.}\label{fig:locw} 
\end{figure}

\begin{figure}[htbp] 
\centerline{\includegraphics[width=.8\textwidth]{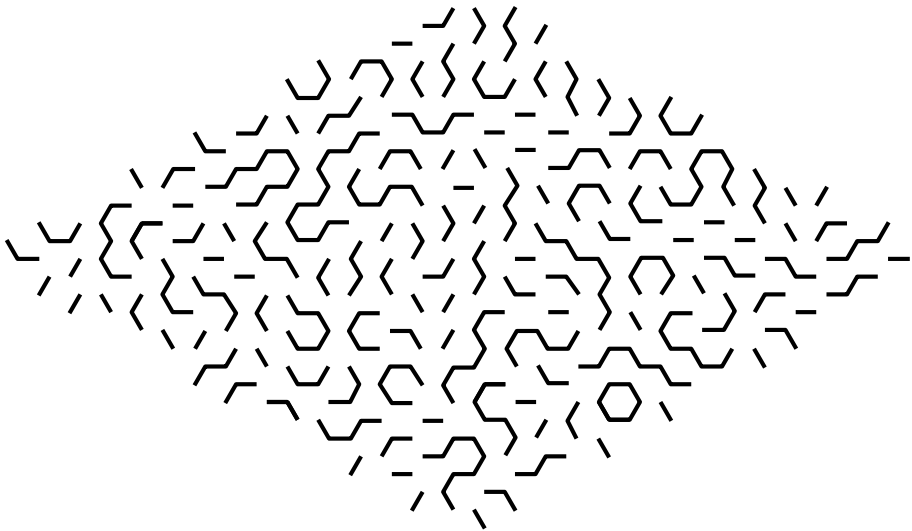}}
 \caption{A realization of the (uniform) \ot model on $\HH_{10}$ with 
 $a=b=c=1$, drawn by MCMC methods. 
 Each vertex has degree either $1$ or $2$, and therefore all components are
 either paths or cycles.}\label{fig:12c} 
 \end{figure}

A vertex $v \in V_n$ is incident to three edges $e_1,e_2,e_3$ of $\HH_n$, which  
are in the respective orientations: horizontal, NW/SE, and NE/SW.
Let $\si \in \Si_n$, and let the \emph{signature} at $v$ be the 
triple $\si|_v:=\si' (e_3)\si'(e_2)\si'(e_1) \in \{0,1\}^3$ considered as a 
word with three letters in the alphabet with two 
letters. Let $a,b,c \in[0,\oo)$ be such that $(a,b,c)\ne (0,0,0)$.
To the signature $\si|_v$ is allocated the weight $w(\si|_v)\in\{0,a,b,c\}$ 
given in Figure \ref{fig:locw}, and the weight function $w:\Si_n\to[0,\oo)$ is defined by
\begin{equation}\label{eq:wt}
w(\si) = \prod_{v\in V_n} w(\si|_v).
\end{equation}
This gives rise to the probability measure $\mu_n$ given by
\begin{equation}\label{eq:pm}
\mu_n(\si) = \frac 1{Z_n} w(\si),
\end{equation}
where 
\begin{equation}\label{eq:pm2}
Z_n=Z_n(a,b,c)= \sum_{\si\in\Si_n}w(\si)
\end{equation}
is the \emph{partition function}. A sample drawn (approximately) from $\mu_{10}$ is 
depicted in Figure \ref{fig:12c}.

It turns out that the weak limit of the sequence $(\mu_n: n \ge 1)$ exists. However, no simple correlation
inequality is known, and the proof of existence of the limit
follows a different route using the relationship
to dimer configurations outlined in Section \ref{ssec:dim}.

\begin{theorem}[{\cite[Thm 6.2]{GL6}}]\label{cl}
The weak limit 
$$
\mo:= \lim_{n\to\oo} \mu_n
$$
exists and is translation invariant.
\end{theorem}

The infinite-volume limit $\mo$ is of course a Gibbs state in the sense that
it satisfies the relevant Dobrushin--Lanford--Ruelle (DLR)
condition (see \cite[Sect.\ 4.4]{G-RCM}). On the other hand,
the structure of the space of such Gibbs measures is unknown.
Neither is it known for which parameter-values $\mo$ is ergodic
(it is not ergodic under the conditions of \cite[Thm 4.9]{ZL2} and Theorem \ref{thm:homo}(b)
of the current paper).

\begin{remark}\label{rem:1}
For edges $e$, $f$ of $\HH$ and for sufficiently large $n$, 
we write $\langle \si_e\si_f\rangle_n$ for the
\emph{two-edge correlation function} of the measure $\mu_n$.
By Theorem \ref{cl}, the limit 
$$
\langle\si_e\si_f\rangle=\lim_{n\to\oo} \langle\si_e\si_f\rangle_n
$$
exists.  In Section \ref{ssec:pt}, we shall consider the limit
$\lim_{|e-f|\to\oo}\langle\si_e\si_f\rangle$ as an order
parameter that is indicative of phase transition. 
\end{remark}

\begin{figure}[htbp] 
\centering
\includegraphics[width=0.9\textwidth]{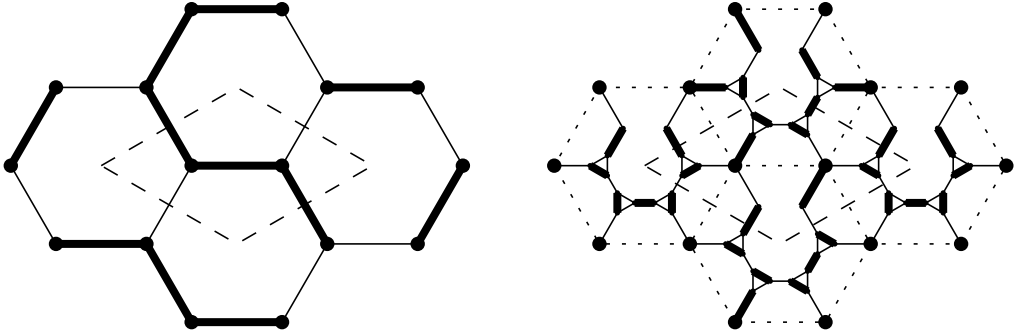} 
\caption{Part of the hexagonal lattice $\HH$ is drawn on the left, together with a \ot
configuration. The graph $\HH_\De$ on the right is obtained by
replacing each face of $\HH$ by a certain `gadget'. 
The left-hand \ot configuration gives rise to a dimer configuration
on the decorated graph on the right, as  described in the text. The fundamental domain 
is outlined as a lozenge, and expanded in Figure \ref{fig:ofd0}.}\label{fig:12con}
\end{figure}

\section{The dimer, Ising, and polygon models}\label{sec:conn}

The relationship between the \ot model  and the dimer model is pivotal to
the study of the former. Dimers are relevant also to the study of Ising models in two dimensions
(see, for example, \cite{ZL1}), and the Ising model gives rise in turn to
a `high temperature' polygon model (see for example,
\cite[p.\ 75]{Bax} and \cite{GJ09,MW73,vdW}). Therefore, the \ot model is connected firmly
to the Ising and polygon models. These connections play roles
in the theory of the \ot model, and are summarized in this section.

\subsection{The dimer model}\label{ssec:dim}

Let $\HH=(\VV,\EE)$ be the hexagonal lattice, and let
$\HH_{\Delta}=(V_{\Delta},E_{\Delta})$ be the \lq decorated graph'
drawn on the right side of Figure \ref{fig:12con}.
The graph in the figure is obtained by
replacing each face $F$ of $\HH$ by a certain `gadget' comprising a path
which is joined to the vertices of $F$ in the manner drawn in the figure. 

A \ot configuration of the left side of Figure \ref{fig:12con} gives rise to a dimer configuration
on the decorated graph on the right, as follows.
Note first that $\VV\subseteq V_{\Delta}$. Each
$v\in \VV$ has three incident edges in $\HH_{\Delta}$, and these 
edges in $\HH_{\Delta}$ may be regarded as the \emph{bisector edges} of the three angles
in $\HH$ at $v$. 

Let $\si\in\Si$. An edge $e\in E_{\Delta}$, incident to a vertex $v\in
\VV$, is designated \emph{present} if and only if the two sides of the corresponding angle
of $\HH$ have the same states, that is, either both are present or both are absent.
Once we have determined the states of the bisector edges of $\HH_{\Delta}$, 
there is a unique extension to a dimer configuration on $\HH_{\Delta}$. See Figure
\ref{fig:12con}. Note that the two \ot configurations $\si,-\si\in\Si$
give rise to the same dimer configuration, and thus the above  correspondence is two-to-one.

Consider now the toroidal graph $\HH_n$ and the corresponding decorated 
graph $\HnD$. The edges of $\HnD$ are weighted, with the edge $e=\langle i,j\rangle$ 
having weight 
\begin{equation}\label{eq:weights}
w_{i,j} = 
\begin{cases} a &\text{if $e$ is a horizontal bisector edge,}\\
 b &\text{if $e$ is a NW/SE bisector edge,}\\ 
 c &\text{if $e$ is a NE/SW bisector edge,}\\
1 &\text{otherwise}.
\end{cases}
\end{equation}
The weight of a dimer configuration is defined as the product of the weights of the edges
that are present, and 
this gives rise (as in \eqref{eq:pm})
to a probability measure $\de_{n,\De}$ on dimer configurations.
Now, $\de_{n,\De}$-probabilities may be represented as weighted counts
of dimer configurations, and such quantities on planar graphs may be computed
by the Pfaffian method of Kasteleyn \cite{Kast61}, Fisher \cite{F61},
and Temperley and Fisher \cite{TF61}. 
This leads to the following limit theorem.

\begin{theorem}[{\cite[Prop.\ 3.3]{ZL2}}]\label{thm:dimer}
Let $a,b,c>0$.  The limit measure 
$$
\de_\De:=\lim_{n\to\oo} \de_{n,\De}
$$ 
exists and is 
translation-invariant and ergodic.
\end{theorem}

Theorem \ref{cl} follows by Theorem \ref{thm:dimer}
and the above correspondence 
between \ot model configurations and dimer configurations. 
Note that the weak limit for \ot measures need not be ergodic.

\begin{figure}[hbtp]
\centerline{
\includegraphics[width=0.56\textwidth]{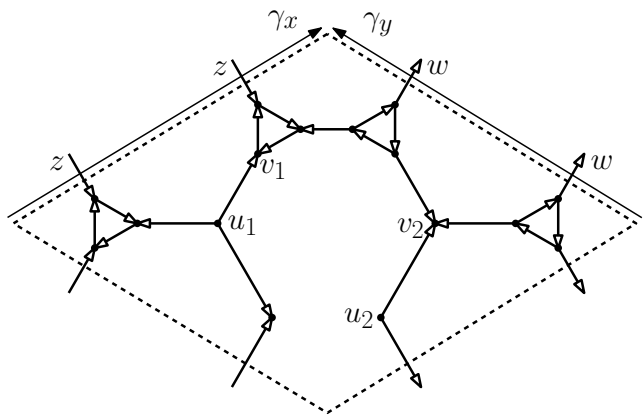}}
\caption{A single fundamental domain of the decorated graph $\HnD$
obtained from the central lozenge of Figure \ref{fig:12con}.
This illustrates the relationship between the fundamental
domain and the original hexagonal lattice $\HH$. Note  the homology cycles
$\g_x$, $\g_y$ of the torus.}
\label{fd}\label{fig:ofd0}
\end{figure}

A $\de_{n,\De}$-probability may be expressed in terms of weighted counts of dimers,
and these are studied
via the Pfaffian representation.
As explained in \cite{GL6} and the references therein, the asymptotics (as $n\to\oo$)
of such Pfaffians depend on the so-called \emph{characteristic polynomial} 
of the model. We do not define the characteristic polynomial here beyond saying 
that it is the determinant $P(z,w)$ of the weighted adjacency matrix of the fundamental domain
of $\HH_\De$, oriented in a `clockwise odd' manner, and 
illustrated in Figure \ref{fig:ofd0}. It is a function of the parameters
$a$, $b$, $c$, and  of two complex variables $z$, $w$. It is shown in
\cite{ZL2} that 
\begin{align*} \label{pzw}
P(z,w) &=
a^4+b^4+c^4+6a^2b^2+6a^2c^2+6b^2c^2-2ab\left(z+\frac{1}{z}\right)\left(a^2+b^2-c^2\right)\\
 &\qq-2ac\left(w+\frac{1}{w}\right)\left(a^2+c^2-b^2\right)-2bc\left(\frac{z}{w}+\frac{w}{z}\right)
 \left(b^2+c^2-a^2\right).\nonumber
 \end{align*}

The \emph{spectral curve} is the zero locus of the characteristic polynomial, that
is, the set of roots of the equation
$P(z,w)=0$. As explained in \cite{ZL2}, it is important to understand the
intersection of the spectral curve with the unit torus 
$$
\TT^2=\bigl\{(z,w)\in \CC^2: |z|=|w|=1\bigr\}.
$$
It turns out that the intersection 
is either empty or is a single real 
point $(1,1)$. Moreover, when $P(1,1)=0$, the zero $(1,1)$ has multiplicity $2$.
Evidently,
\begin{equation}\label{eq:303}
P(1,1) = (a^2+b^2+c^2-2ab-2bc-2ac)^2,
\end{equation}
and therefore the spectral curve intersects $\TT^2$ if and only if
\begin{equation}\label{eq:spcurve}
\sqrt a \pm \sqrt b \pm \sqrt c=0.
\end{equation}
We shall return to this equation in the study of phase
transition in Section \ref{ssec:pt}.

\subsection{The half-edge graph}\label{ssec:heg}

When considering correlations, it will be convenient to work on a graph derived from 
the hexagonal lattice by replacing each edge by two half-edges.
Let $A \HH_n=(A V_n, A E_n)$ be the graph derived 
from $\HH_n=(V_n,E_n)$ by adding
a vertex at the midpoint of each edge in $E_n$. Let $M E_n=\{M e: e \in E_n\}$ be the set
of such midpoints, and $A V_n = V_n \cup M E_n$.  The edges $A E_n$
are the half-edges of $E_n$, each being of the form $\langle v, Me\rangle$
for some $v \in V_n$ and incident edge $e \in E_n$.

The \ot model on $\HH_n$ can be viewed as a spin-model on the set $ME_n$
of midpoints, as we explain next.

\subsection{The Ising model}\label{ssec:ising}

It turns out that the \ot model is the marginal of a
certain Ising-type model on
the half-edge graph $A\HH_n$ of Section \ref{ssec:heg}, that 
is reminiscent of the Edwards--Sokal coupling
of the Potts and random-cluster measures (see \cite[Sect.\ 1.4]{G-RCM}).
It is constructed via a weight function on configuration space,
using weights that are permitted in general to be $\CC$-valued. 

Let $\Si^\re=\{-1,+1\}^{ME_n}$ and $\Si^\rv =\{-1,+1\}^{V_n}$.
An edge $e\in E_n$ is identified with the element of 
$ME_n$ at its centre.
A spin-vector is a pair $(\si^\re,\si^\rv)\in \Si^\re \times \Si^\rv$ with 
$\si^\re=(\si_{v,s}: v\in V_n,\ s=a,b,c)$ and
$\si^\rv=(\si_v: v \in V_n)$, where
$\si_{v,a}$, $\si_{v,b}$, $\si_{v,c}$ denote the spins on midpoints of the corresponding 
edges of given types incident to $v\in V_n$ (see Figure \ref{fig:hex0}).
We allocate the (possibly complex)  weights
\begin{equation}\label{eq:ss}
w(\si^\rv,\si^\re) = \prod_{v\in V_n}(1+\eps_a\si_v\si_{v,a})(1+\eps_b\si_v\si_{v,b})(1+\eps_c\si_v\si_{v,c}),
\end{equation}
where $\eps_a, \eps_b, \eps_c\in\CC$ are 
constants associated with horizontal, NW/SE, and NE/SW edges, respectively. 
In \eqref{eq:ss}, each factor $1+\eps_s \si_v\si_{v,s}$ ($s=a,b,c$) 
corresponds to a half-edge of $\HH _n$.

When considering the relationship to the \ot model, it
will be convenient to choose the constants $\eps_s$, $s=a,b,c$, as follows. Let $a,b,c>0$, and
\begin{equation}
A=\frac{a-b-c}{a+b+c},\q 
B=\frac{b-a-c}{a+b+c},\q 
C=\frac{c-a-b}{a+b+c},\label{abc}
\end{equation} 
where we assume  for simplicity that $ABC\ne 0$.
The appropriate values of the $\eps_s$ are
\begin{equation}\label{eq:half2-}
\eps_a=\sqrt{\frac{BC}{A}},\q \eps_b=\sqrt{\frac{AC}{B}},\q \eps_c=\sqrt{\frac{AB}{C}}.
\end{equation}

\begin{theorem}[{\cite[Sect.\ 4]{GL6}}]\mbox{}
\begin{letlist}

\item \emph{Marginal on $V_n$}.
For  $\si^\rv\in\Si^\rv$, $w(\si^\rv):=
\sum_{\si^\re\in\Si^\re} w(\si^\rv,\si^\re)$ satisfies
\begin{equation}\label{eq:pfss2}
 w(\si^\rv) = 2^{|E_n|}
 \prod_{g=\langle u,v\rangle\in E_n}
\left(1+\eps_g^2\si_u\si_v\right).
\end{equation}
That is, the marginal weights on $\Si^\rv$ are those
of an Ising-type model on $\HH_n$ 
with (possibly complex) edge-interactions. Here, $\eps_g$ denotes the parameter associated with 
edge $g=\langle u,v\rangle$.

\item \emph{Marginal on $ME_n$}. Let $a,b,c>0$ and assume the
$\eps_s$ satisfy \eqref{abc}--\eqref{eq:half2-} where $ABC\ne 0$.
For $\si^\re\in\Si^\re$, $w(\si^\re):=
\sum_{\si^\rv\in\Si^\rv} w(\si^\rv,\si^\re)$ satisfies
\begin{align*}
w(\si^\re) &=\prod_{v\in V_n}
\bigl(1+A\si_{v,b}\si_{v,c}+B\si_{v,a}\si_{v,c}+
C\si_{v,a}\si_{v,b}\bigr)\\
&\propto \mu_n(\si^\re).
\end{align*}
That is, the marginal weights on $\Si^\re$ are proportional to those of the \ot model.

\item \emph{Two-edge correlation function}.
Let $e=\langle u,v\rangle, f=\langle x,y\rangle \in E_n$.
Subject to the notation of part {\rm(b)} above, the two-edge correlation
function of Remark \ref{rem:1} satisfies
\begin{equation}\label{eq:corr4}
\langle\si_e\si_f\rangle_n =  
\left.\sum_{\si^\rv\in\Si^\rv} D_{e,f}(\si^\rv) w(\si^\rv)\right/ \sum_{\si^\rv\in\Si^\rv}w(\si^\rv),
\end{equation}
where
$$
D_{e,f}(\si^\rv) = 
\frac{\eps_e(\si_u+\si_v)\eps_f(\si_x+\si_y)}
{(1+\eps_e^2)(1+\eps_f^2)}.
$$ 
\end{letlist}
\end{theorem}

If the weights $w(\si^\rv)$ of \eqref{eq:pfss2}
are real and non-negative (which they are not in general), 
the ratio on the right side of \eqref{eq:corr4} may be interpreted as an expectation. 
The weights
$w(\si^\rv)$ correspond to a ferromagnetic Ising model if and only
the edge-weights of \eqref{eq:pfss2} satisfy $0<\eps_g^2<1$.
Unfortunately, this  never occurs with the $\eps_g$ 
derived from the \ot model as in \eqref{abc}--\eqref{eq:half2-}.
If, however, one assumes that the $a$, $b$, $c$
satisfy the \emph{acute angle condition}
$$
a^2<b^2+c^2,\q b^2< c^2+a^2, \q c^2<a^2+b^2,
$$
then $-1<\eps_g^2<0$, and the corresponding Ising model is antiferromagnetic.
Since $\HH$ is bipartite, this process may be transformed into a ferromagnetic system
by changing the sign of every other vertex
(see \cite[p.\ 17]{GHM}), and this transformation greatly facilitates its analysis.

The above Ising model may be regarded as
a special case of the eight-vertex model of Lin and Wu \cite{LW90}.

\subsection{The hexagonal polygon model}

Let $\HH_n=(V_n,E_n)$ as before, and let $\Pi_n=\{0,1\}^{E_n}$.
The sample space of the polygon model on $\HH_n$ is the subset $\Pip_n\subseteq\Pi_n$ 
containing all $\pi=(\pi(e): e \in E_n) \in\Pi_n$ such that 
\begin{equation}\label{eq:polycond}
\sum_{e\ni v} \pi(e)\q \text{ is either $0$ or $2$},\qq  v \in V_n.
\end{equation} 
Each $\pi\in\Pip_n$ may be considered as a union of vertex-disjoint cycles of $\HH_n$, together
with isolated vertices. We identify $\pi\in\Pi_n$ with the set $\{e\in E_n:\pi(e)=1\}$  of `open' edges
under $\pi$.  Thus \eqref{eq:polycond} requires that every vertex is incident
to an even number of open edges.

Let $\eps_a, \eps_b, \eps_c\in \CC$. To the configuration $\pi\in\Pip_n$, we assign the 
(possibly complex) weight
\begin{equation}\label{eq:polywt}
w(\pi)= \eps_a^{2|\pi(a)|}\eps_b^{2|\pi(b)|} \eps_c^{2|\pi(c)|},
\end{equation}
where $\pi(s)$ is  the set of open $s$-type edges. 
The weight function $w$ gives rise to the partition function
\begin{equation*}\label{eq:pmp}
\Znp=\sum_{\pi\in\Pip_n} w(\pi).
\end{equation*}

Let $a,b,c>0$.
We now choose the constants $\eps_s$, $s=a,b,c$, according to 
\eqref{abc}--\eqref{eq:half2-}, where it is assumed that $ABC\ne 0$.
The corresponding polygon model is related to the high-temperature expansion
of the Ising-type model of Section \ref{ssec:ising} (see, for example, 
\cite[eqn 5.1]{G-FF} and \cite[Thm 1.7]{GJ09}).

In considering correlation functions, it is convenient  to view the polygon model  on 
the half-edge graph $A\HH_n$ of Section \ref{ssec:heg}.
A polygon configuration on
$\HH_n$ induces a polygon configuration on $A \HH_n$, namely a 
subset of $AE_n$ with the property that every vertex in $AV_n$ has even degree. 
For an $a$-type edge $e\in E_n$,
the two half-edges of $e$ have weight $\eps_a$ each (and similarly for $b$- and $c$-type edges).
The weight function $w$ of \eqref{eq:polywt} may now be expressed as
$$
w(\pi)= \eps_a^{|\pi(a)|}\eps_b^{|\pi(b)|} \eps_c^{|\pi(c)|},
$$
where $\pi(s)$ is the set of open half-edges of type $s$, 
as $\pi$ ranges over polygon configurations on $A\HH_n$.

Let $e, f\in ME_n$ be distinct midpoints of $A\HH_n$, and let $\Pipi$ be the subset 
of all $\pi\in\{0,1\}^{AE_n}$
such that: (i) every $v\in AV_n$ with $v\ne e,f$ is incident to an even number of open half-edges,
and (ii)  the midpoints $e$ and $f$ are incident to exactly one open half-edge.
Let
\begin{equation}\label{eq:edgecor2}
M_n(e,f)=\frac{Z_{n,e\lra f}}{\Znp},
\end{equation}
where 
\begin{equation}\label{eq:edgecor3}
Z_{n,e\lra f} := \sum_{\pi\in\Pipi} \eps_a^{|\pi(a)|}\eps_b^{|\pi(b)|} \eps_c^{|\pi(c)|}.
\end{equation}

\begin{theorem}[\cite{GL7}]
Subject to \eqref{eq:half2-}, the two-edge correlation 
function $\langle \si_e\si_f\rangle_n$ of the \ot model 
on $\HH_n$ satisfies $\langle\si_e\si_f\rangle_n=M_n(e,f)$.
\end{theorem}

The polygon model with general parameters is studied in \cite{GL7}.

\section{Phase transition}\label{sec:pt}

We discuss two forms of phase transition for the \ot model. Of these, the first
concerns the existence (or not) of infinite `homogeneous' clusters of $\HH$ containing 
vertices of the same type, and the second considers as order parameter 
the limiting two-edge correlation function. 
Thus, the first studies the geometry of the model,
and the second its correlation structure.
There may exist other forms of phase transition,
as yet unstudied.

\subsection{Occurrence of paths}\label{ssec:occ}

Every connected component (`cluster') in a realization of the \ot model
is either a self-avoiding path or a cycle (see Figure \ref{fig:12c}). 
It turns out that all such clusters are $\mo$-a.s.\ finite,
when $a,b,c>0$. We formalize this statement in this subsection, and begin with an exact formula.
Using the correspondence between
\ot model configurations on $\HH$ and dimer configurations on $\HH_{\Delta}$,
as described in Section \ref{ssec:dim}, we have the following.

\begin{theorem}[{\cite[Thm 3.4]{ZL2}}]\label{thm:41}
Let $a,b,c>0$ and let $\mo$ be the limit \ot measure of Theorem \ref{cl}. Let
$\pa$ be a self-avoiding path of $\HH_n$ containing $l+1$ edges, and write 
$E_{\pa}=\{e_k=\langle u_k,v_k\rangle: 1\le k \le l\}$ for the set of 
bisector edges of $\HH_\De$ encountered along $\pa$, as in Figure \ref{fig:12con}.
Then 
\begin{equation*}
\mo(\text{\rm $\pa$ is present})
=\tfrac{1}{2}\left(\prod_{k=1}^{l}w_{e_k}\right)\left|\Pf K^{-1}_\pa\right|,
\end{equation*} 
where $w_{e}$ is the weight of the edge $e$ in
$\HH_{\Delta}$, $K^{-1}_\pa$ is the submatrix of the inverse of
the weighted adjacency matrix of $\HH_{\Delta}$ with rows and columns indexed by
$u_1,v_1,\dots,u_l,v_l$, and $\Pf M$ is the Pfaffian of the matrix $M$.
\end{theorem}

The mass transport principle, introduced in \cite{blps,hagg97}, is a valuable tool in
the study of interacting systems including percolation 
and self-avoiding walks on Cayley graphs, see \cite{blps1,GL-Cayley,hp}. 
It may also be applied to the \ot model, where it is
used to prove the following.

\begin{theorem}[{\cite[Thm 2.4]{ZL2}}]\label{thm:no}
If $a,b,c>0$, we have that
$$
\mo\bigl(\text{\rm there exists an infinite path}\bigr)=0.
$$
\end{theorem}

\subsection{Existence of infinite homogeneous clusters}\label{ssec:inf}

The concept of `phase transition' hinges on the non-smoothness of some so-called
`order parameter'.  For the Ising model, one may take as order
parameter the magnetization at 
the origin in the infinite-volume measure with $+$ boundary conditions.
This corresponds in the universe of percolation and the random-cluster model 
(see \cite{G-RCM}) to
studying whether or not there there exists an infinite cluster.
By Theorem \ref{thm:no}, the \ot model possesses no infinite cluster for any values
of $a,b,c>0$. `Clusters' may however be defined in another manner.

Let $\si\in\Si$. Each vertex of
$\HH$ has a random type, given in Figure \ref{fig:locw}. 
For $s\in\{a,b,c\}$, a \emph{type-$s$} cluster is a maximal
connected subgraph of $\HH$ every vertex of which has type $s$. 
As illustrated in the figure, type $s$ comes in two sub-types; for example, 
a type-$a$ vertex has signature either $001$ or $110$. Thus, one
may speak of a $001$-cluster, etc. By examining the figure again,
it is seen that a type-$a$ cluster is either a $001$-cluster
or a $110$-cluster, but may not contain vertices of both
types (and similarly for types $b$ and $c$).
A \emph{homogeneous cluster} of $\si\in\Si$ is a $w$-cluster for some
$w\in\{0,1\}^3$, $w\ne 000,111$. We concentrate now
on the existence (or not) of an infinite homogeneous cluster.

\begin{theorem}\label{thm:homo}
Let $a,b,c>0$.
\begin{letlist}
\item {\rm(\cite[Thm 1.1]{ZL3})} 
Let $w\in\{0,1\}^3$, $w\ne 000,111$. The number $I_w$ of infinite
$w$-clusters is $\mo$-a.s.\ no greater than $1$.  

\item {\rm(\cite[Thm 4.4, Prop.\ 4.7]{ZL2})} 
Fix $b,c>0$. For sufficiently small $a>0$, there exists $\mo$-a.s.\ no
infinite type-$a$ cluster. 
For sufficiently large $a$,
the $\mo$-probability that the origin belongs to an infinite
type-$a$ cluster is strictly positive.

\end{letlist}
\end{theorem}

Part (a) is proved in \cite{ZL3} using an adaptation of the method of  Burton and Keane 
(see \cite{bk89,ns81}) to the \ot model, 
subject to the complication that $\mo$ does not have the so-called 
`finite energy property'. It is unknown whether infinite
$w$-clusters and $w'$-clusters can coexist with $w \ne w'$.

Part (b) indicates the existence of a phase transition. It is not known if
there exists a single critical point $\ac=\ac(b,c)$ for the given property. Furthermore, there
is currently no indication of the exact value of such a point. 

\subsection{Non-analyticity of the two-edge correlation function}\label{ssec:pt}

\begin{figure}[htbp]
  \centering
\includegraphics*[scale=1]{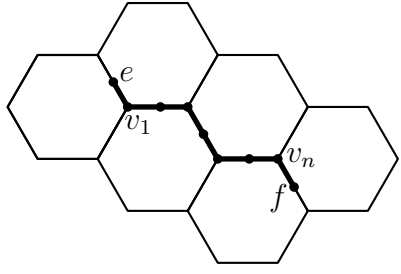}
   \caption{A path $\pa$ comprising horizontal and NW/SE mid-edges, 
   connecting the midpoints of two NW/SE edges $e$ and $f$.}\label{fig:lef2}
\end{figure}

For distinct edges $e,f\in \EE$, we write 
$$
\langle \si_e\si_f\rangle = \lim_{n\to\oo} \langle \si_e\si_f\rangle_n,
$$ 
which exists by Theorem \ref{cl}, (see \cite[Thm 6.2]{GL6}).
We consider here the asymptotic behaviour of $\langle \si_e\si_f\rangle$
as $|e-f|\to\oo$. The behaviour of this limit is unknown in general, but
a great deal is known if $e$ and $f$ are related in the `diagonal' manner
of the forthcoming assumption \eqref{eq:condition0}, as illustrated in
Figure \ref{fig:lef2}.

\begin{theorem}[{\cite[Thm 3.1]{GL6}}]\label{thm:main}
Let $a,b, c >0$, and let $e,f \in \EE$ be NW/SE edges such that: 
\begin{equation}\label{eq:condition0}
\begin{aligned}
&\text{there exists a path $\pa=\pa(e,f)$ of $\HH_n$ from $e$ to $f$}\\
& \text{using only horizontal and NW/SE half-edges}.
\end{aligned}
\end{equation}
\begin{letlist}
\item 
Let $a \ge b>0$. For $c>0$ satisfying 
$$
\text{\rm either $\sqrt a > \sqrt b + \sqrt c$\q or\q $\sqrt{c}>\sqrt{a}+\sqrt{b}$}
$$
except possibly on some set of isolated points, the limit
$\lim_{|e-f|\to\oo}\langle \si_e\si_f\rangle^2$ exists and is non-zero.
\item If $a \ge b>0$ and 
$$
\sqrt a-\sqrt b < \sqrt c<\sqrt a+\sqrt b,
$$
then $\langle \si_e\si_f\rangle \to 0$ as $|e-f|\to\oo$.
\end{letlist}
\end{theorem}

By Theorem \ref{thm:main}, when $a\ge b\ge c>0$, the phase transition occurs when 
$\sqrt{a}=\sqrt{b}+\sqrt{c}$. The proof is along the following lines.
The square 
$\langle \si_e\si_f \rangle^2$ 
of the two-edge correlation may be expressed
as the determinant of an explicit block Toeplitz matrix
(this is where \eqref{eq:condition0} is used); see the forthcoming Lemma \ref{lem:pf2}.
Its limit as
$|e-f|\to\oo$ is given by Widom's theorem (see \cite{HW0,HW} 
and \cite[Thm 8.7]{GL6}) as the determinant of the 
limiting (infinite) block
Toeplitz matrix. This determinant  is complex analytic with respect to the parameters $a$, $b$, $c$ 
except when the
spectral curve has a unique real zero on the unit torus. As remarked at the end of 
Section \ref{ssec:dim}, the last occurs if and only if $\sqrt a \pm \sqrt b \pm \sqrt c=0$.
When $a \ge b \ge c$, this equation becomes $\sqrt a-\sqrt b -\sqrt c=0$. 

The `isolated points' of part (b) arise through the use in the proof of the fact that, 
for an analytic function $\La$, either $\La\equiv 0$ or the zeros of $\La$ are isolated.

It is easily seen from Figures \ref{fig:locw} and \ref{fig:lef2} that
\begin{equation*}\label{eq:extreme}
\langle \si_e\si_f\rangle = 1\qq\text{if either}\q a,b>0,\ c=0,\q \text{or}\q a=b= 0,\ c>0,
\end{equation*}
and part (b) of Theorem \ref{thm:main} follows by the analyticity.
For part (a), one uses the representation of the \ot model as the Ising model of
Section \ref{ssec:ising}. 

The key step in the proof of Theorem \ref{thm:main} is the following exact formula.
Let $Y_1$ be the $2\times 2$ matrix
\begin{equation*}
Y_1=\begin{pmatrix} 0&1\\-1&0 \end{pmatrix},
\end{equation*}
and let $Y_{2k}$ be the $4k\times 4k$ block diagonal matrix with 
diagonal $2\times 2$ blocks equal to $Y_1$. That is, 
\begin{equation*}
Y_{2k}=\begin{pmatrix}
Y_1&0&\cdots&0\\
0&Y_1&\cdots&0\\
\vdots&\vdots&\ddots&\vdots\\
0&0&\cdots&Y_1
\end{pmatrix}.
\end{equation*}

\begin{lemma}[{\cite[Lemma 8.4]{GL6}}]\label{lem:pf2}
Suppose the path $\pa=\pa(e,f)$ of \eqref{eq:condition0} passes $2k$ bisector edges of $\HH_\De$.
We have that
\begin{equation*}\label{eq:236}
\langle \si_e\sigma_f \rangle
=\Pf[Y_{2k}+2cK_\pa^{-1}],
\end{equation*}
where $K^{-1}_\pa$ is as in Theorem \ref{thm:41}.
\end{lemma}

\section*{Acknowledgements} 
This work was supported in part by the Engineering
and Physical Sciences Research Council under grant EP/103372X/1. ZL acknowledges
support from the Simons Foundation under grant $\#$351813. The authors thank
the anonymous referee for the careful reading and useful comments.

\bibliography{12review6} 
\bibliographystyle{amsplain}
\end{document}